\documentclass[11pt,twoside]{article}
\usepackage[letterpaper,left=1.5in,top=1in,bottom=1in,right=1.5in]{geometry}
\newcommand\Z{\mathbf Z}
\newcommand\R{\mathbf R}
\newcommand\Q{\mathbf Q}
\newcommand\round[2][m]{\lfloor #2\rceil_{#1}}
\newcommand{\st}{\mathrm{s.t.}}
\usepackage{amsmath,amsfonts,amsthm,amsbsy,url,enumerate}
\let\ve=\mathbf  
\newcommand\vealpha{{\boldsymbol{\alpha}}}
\newcommand\vebeta{{\boldsymbol{\beta}}}

\newtheorem{theorem}{Theorem}[section]
\newtheorem{lemma}[theorem]{Lemma}

\newtheorem{corollary}[theorem]{Corollary}
\theoremstyle{definition}

\newtheorem{example}[theorem]{Example}
\newtheorem{algorithm}[theorem]{Algorithm}
\usepackage{hyperref}

\makeatletter
\let\@fnsymbol=\@arabic
\makeatother
\title{FPTAS for mixed-integer polynomial optimization with a fixed
  number of variables}
\author{J.~A.~De Loera\thanks{Address: University of California, Dept.~of Mathematics,
    Davis CA 95616, USA; E-Mail Address: \texttt{deloera@math.ucdavis.edu}.}
  \and R.~Hemmecke\thanks{Address: Otto-von-Guericke-Universit\"at Magdeburg, FMA/IMO,
    Universit\"atsplatz 2, 39106 Magdeburg, Germany; E-Mail Address:
    \texttt{hemmecke@imo.math.uni-magdeburg.de}.}
  \and M.~K\"oppe\thanks{Address: Otto-von-Guericke-Universit\"at Magdeburg, FMA/IMO,
    Universit\"atsplatz 2, 39106 Magdeburg, Germany; E-Mail Address:
    \texttt{mkoeppe@imo.math.uni-magdeburg.de}.}
  \and R.~Weismantel\thanks{Address: Otto-von-Guericke-Universit\"at Magdeburg, FMA/IMO,
    Universit\"atsplatz 2, 39106 Magdeburg, Germany; E-Mail Address:
    \texttt{weismant@imo.math.uni-magdeburg.de}.}}

\date{$\relax$Revision: 1.40 $ - \ $Date: 2005/05/31 10:03:50 $ $}

\begin{document}
\maketitle

\begin{abstract}
  \noindent We show the existence of an FPTAS for the problem of maximizing a
  non-negative polynomial over mixed-integer sets in convex
  polytopes, when the number of variables is fixed.
\end{abstract}

\section{Introduction}
A well-known result by H.W. Lenstra Jr.\ states that {\em linear} mixed
integer programming problems with fixed number of variables can be
solved in polynomial time on the input size \cite{Len83}. It is a natural
question to ask what is the computational complexity, when the number
of variables is fixed, of the {\em non-linear} mixed integer problem
\begin{subequations}
  \label{eq:mixed-problem}
  \begin{align}
    \max\quad & f(x_1,\dots,x_{d_1}, z_1,\dots,z_{d_2}) \\
    \st\quad & A\ve x + B\ve z \leq \ve b \label{eq:mixed-problem-a}\\
    & x_i \in \R && \text{for $i = 1, \dots, d_1$}, \\
    & z_i \in \Z && \text{for $i = 1, \dots, d_2$}, \label{eq:mixed-problem-b}
  \end{align}
\end{subequations}
where $f$ is a polynomial function of maximum total degree~$D$ with
rational coefficients, and $A\in\Z^{m\times d_1}$, $B\in\Z^{m\times
  d_2}$, $\ve b\in\Z^m$ (here we assume that $A\ve x + B\ve z \leq \ve
b$ describes a convex polytope, which we denote by~$P$).

It was well-known that continuous polynomial optimization over
polytopes, without fixed dimension, is NP-hard and that an {\small FPTAS} is
not possible. Indeed the max-cut problem can be modeled as minimizing
a quadratic form over the cube $[-1,1]^d$ \cite{Hastad:inapprox97}.
More strongly, it turns out that, even for dimension two and total
degree of $f$ four, problem (\ref{eq:mixed-problem}) is an NP-hard problem too
\cite{deloera-hemmecke-koeppe-weismantel:intpoly-fixeddim}. Thus the
best we can hope for, even for fixed dimension, is an approximation result. This paper presents
the best possible such result:

\begin{theorem}
  \label{th:mipo-fptas}
  Let the dimension $d=d_1+d_2$ be fixed.  
  \begin{enumerate}[\rm(a)]
  \item There exists a fully
  polynomial time approximation scheme ({\small FPTAS}) for the
  optimization problem~\eqref{eq:mixed-problem} for all polynomial
  functions $f\in\Q[x_1,\dots,x_{d_1},\allowbreak z_1,\dots,z_{d_2}]$
  that are non-negative on the feasible
  region~\textup{(\ref{eq:mixed-problem-a}--\ref{eq:mixed-problem-b})}.
  (We assume the encoding length of~$f$ is at least as large as its
  maximum total degree.)
  
\item Moreover, the restriction to non-negative polynomials is necessary, as there
  does not even exist a polynomial time approximation scheme ({\small PTAS})
  for the maximization of {\em arbitrary} 
  polynomials over mixed-integer sets in polytopes, even for fixed dimension
  $d\geq2$.
\end{enumerate}
\end{theorem}

The proof of \autoref{th:mipo-fptas} is presented in 
\autoref{section:mainthm}.  As we will see, \autoref{th:mipo-fptas} is
a non-trivial consequence of the existence of {\small FPTAS} for the
problem of maximizing a non-negative polynomial with integer
coefficients over the lattice points of a convex rational polytope.
That such {\small FPTAS} indeed exists was recently settled in
our paper~\cite{deloera-hemmecke-koeppe-weismantel:intpoly-fixeddim}.  
The knowledge of
paper~\cite{deloera-hemmecke-koeppe-weismantel:intpoly-fixeddim} is not
necessary to understand this paper but,
for convenience of the reader, we include a short summary in an appendix.
Our
arguments, however, are independent of which {\small FPTAS} is used in the integral case. 
Our results come to complement other approximation schemes
investigated for continuous variables and fixed degree
(see \cite{deklerk-laurent-parillo:ptas-polynomial-simplex} and
references therein).  

One interesting property of our {\small FPTAS} is that it
does not depend on dynamic programming, unlike most known {\small FPTAS} (see
comments and references in the introduction
of~\cite{woeginger-dynamic-programming}). 
Instead our main approach is to use grid refinement in order to approximate the
mixed-integer optimal value via auxiliary pure integer problems.
One of the difficulties on constructing
approximations is the fact that not every sequence of grids whose
widths converge to zero leads to a convergent sequence of optimal
solutions of grid optimization problems.  This difficulty is addressed in 
\autoref{section:qualityofgrid}. 
In \autoref{section:constant} we develop techniques for bounding
differences of polynomial function values.
Finally, \autoref{section:mainthm} contains the proof of
\autoref{th:mipo-fptas}.

\section{Grid approximation results}
\label{section:qualityofgrid}

An important step in the development of an {\small FPTAS} for the mixed-integer
optimization problem is the reduction of the mixed-integer
problem~\eqref{eq:mixed-problem} to an auxiliary optimization problem over a
lattice $\frac1m \Z^{d_1}\times \Z^{d_2}$.  To this end, we consider the
\emph{grid problem} 
\begin{equation}
  \label{eq:grid-problem}
  \begin{aligned}
    \max\quad & f(x_1,\dots,x_{d_1},z_1,\dots,z_{d_2}) \\
    \st\quad & A\ve x + B\ve z\leq \ve b \\
    & x_i \in \tfrac 1m \Z &&\text{for $i=1,\dots,d_1$},\\
    & z_i \in \Z &&\text{for $i=1,\dots,d_2$}.\\
  \end{aligned}
\end{equation}
We can solve this problem approximately using the integer {\small FPTAS}
(\autoref{lemma:bounds}):
\begin{corollary}\label{cor:gridproblem-approximation}
  For fixed dimension $d=d_1+d_2$ there 
  exists an algorithm with running time polynomial in $\log m$, in the
  encoding length of~$f$ and of~$P$, in the maximum total degree~$D$ of~$f$, and
  in~$\frac1\epsilon$ for computing a 
  feasible solution $(\ve x^m_\epsilon, \ve z^m_\epsilon)\in
  P\cap\bigl(\frac1m\Z^{d_1}\times \Z^{d_2}\bigr)$ to the grid
  problem~\eqref{eq:grid-problem}, where $f$~is non-negative on the feasible
  region, with
  \begin{equation}\label{eq:gridproblem-approximation}
    f(\ve x^m_\epsilon, \ve z^m_\epsilon) \geq (1-\epsilon) f(\ve x^m,\ve z^m),
  \end{equation}
  where $(\ve x^m, \ve z^m)\in
  P\cap\bigl(\frac1m\Z^{d_1}\times \Z^{d_2}\bigr)$ is an optimal solution
  to~\eqref{eq:grid-problem}.
\end{corollary}
\begin{proof} 
  We apply \autoref{lemma:bounds} to the pure
  integer optimization problem:  
  \begin{equation} 
    \label{eq:pure-integer-problem}
    \begin{aligned}
      \max\quad & \tilde f(\ve{\tilde x}, \ve z) \\
      \st\quad  & A\ve{\tilde x} + m B\ve z \leq m\ve b \\
      & \tilde x_i \in \Z &&\text{for $i=1,\dots,d_1$}, \\
      & z_i \in \Z &&\text{for $i=1,\dots,d_2$},
    \end{aligned} 
  \end{equation}
  where $\tilde f(\ve{\tilde x},\ve z) := m^D f( \tfrac1m \ve{\tilde x}, \ve
  z)$ is a polynomial function with integer coefficients. 
  Clearly the binary encoding length of the coefficients 
  of~$\tilde f$ increases by at most $\lceil D\log m \rceil$, compared to the
  coefficients 
  of~$f$.  Likewise, the encoding length of the coefficients
  of~$mB$ and~$m\ve b$ increases by at most $\lceil \log m\rceil$.
  By Theorem 1.1
  of~\cite{deloera-hemmecke-koeppe-weismantel:intpoly-fixeddim}, there exists
  an algorithm with running time polynomial in the encoding length of~$\tilde
  f$ and of~$A\ve x + m B\ve z \leq m \ve b$, the maximum total degree~$D$,
  and~$\frac1\epsilon$ for computing a feasible solution $(\ve x^m_\epsilon, \ve z^m_\epsilon)\in
  P\cap\bigl(\frac1m\Z^{d_1}\times \Z^{d_2}\bigr)$ such that $\tilde f(\ve
  x^m_\epsilon, \ve z^m_\epsilon) \geq (1-\epsilon) \tilde f(\ve x^m,\ve z^m)$,
  which implies~\eqref{eq:gridproblem-approximation}.
\end{proof} 

One might be tempted to think that for large-enough choice of~$m$, we
immediately obtain an approximation to the mixed-integer optimum with
arbitrary precision.  However, this is not true, as the following example
demonstrates. 

\begin{example} 
  Consider the mixed-integer optimization problem 
  \begin{equation} 
    \label{eq:example-slice-not-fulldim}
    \begin{aligned}
      \max\quad& 2z - x\\
      \st\quad& z \leq 2x \\
      & z \leq 2(1-x) \\
      & x\in\R_+,\ z\in\{0,1\}, 
    \end{aligned} 
  \end{equation}
  whose feasible region consists of the point $(\frac12,1)$ and the segment
  $\{\,(x,0): x\in[0,1]\}$.  The unique optimal solution
  to~\eqref{eq:example-slice-not-fulldim} is $x=\frac12$, $z=1$.  Now consider
  the sequence of grid approximations of~\eqref{eq:example-slice-not-fulldim}
  where $x\in \frac1m \Z_+$.  For even $m$, the unique optimal solution to the
  grid approximation is $x=\frac12$, $z=1$.  However, for odd $m$, the unique
  optimal solution is $x=0$, $z=0$.  Thus the full sequence of the optimal
  solutions to the grid approximations does not converge since it has two limit
  points. 
\end{example} 

However we can prove that it is possible to construct, in polynomial time, a
subsequence of finer and finer grids that contain a lattice point $(\ve x,\ve
z^*)$ that is arbitrarily close to the mixed-integer optimum $(\ve x^*,\ve
z^*)$.  This is the central statement of this section and a basic building
block of the approximation result.
\begin{theorem}[Grid Approximation]
  \label{lemma:grid-approximation-lowerdim}
  Let $d_1$ be fixed.
  Let $P = \{\, (\ve x,\ve z)\in\R^{d_1+d_2} : A\ve x + B\ve z \leq \ve
  b\,\}$, where $A\in\Z^{m\times d_1}$, $B\in\Z^{m\times d_2}$.  Let $M\in\R$
  be given such
  that $P \subseteq \{\, (\ve x,\ve z)\in\R^{d_1+d_2} : |x_i| \leq M \text{ for
  } i=1,\dots,d_1 \,\}$.  
  There exists a polynomial-time algorithm to compute a number $\Delta$
  such that for every $(\ve x^*,\ve z^*)\in P\cap(\R^{d_1}\times \Z^{d_2})$
  and $\delta>0$ the following property holds:
  \begin{quote} 
    Every lattice $\frac{1}{m}\Z^{d_1}\times \Z^{d_2}$ for $m=k\Delta$ and
    $k\geq\frac1\delta (d_1+1) M$ contains a lattice point $(\ve x,\ve z^*)\in
    P\cap \bigl(\frac{1}{m}\Z^{d_1}\times \Z^{d_2}\bigr) $ with $\|\ve x-\ve
    x^*\| \leq\delta$. 
  \end{quote} 
\end{theorem}
\autoref{lemma:grid-approximation-lowerdim} follows directly from the
next two lemmas.

\begin{lemma}[Integral Scaling Lemma]
  \label{lemma:integral-scaling}
  Let $P = \{\, (\ve x,\ve z)\in\R^{d_1+d_2} : A\ve x + B\ve z \leq \ve
  b\,\}$, where $A\in\Z^{m\times d_1}$, $B\in\Z^{m\times d_2}$.
  For fixed $d_1$, there exists a polynomial time algorithm to compute a number
  $\Delta\in\Z_{>0}$ such that for every $\ve z\in\Z^{d_2}$ the polyhedron 
  \begin{displaymath}
    \Delta P_{\ve z} = \bigl\{\, \Delta\ve x : (\ve x,\ve z) \in P\,\bigr\}
  \end{displaymath}
  is integral.  In particular, the number $\Delta$ has an encoding
  length that is bounded by a polynomial in the encoding length
  of~$P$.
\end{lemma}
\begin{proof}
  Because the dimension $d_1$ is fixed, there exist only polynomially many
  simplex bases of the system $A \ve x \leq \ve b - B\ve z $, and they can be
  enumerated in polynomial time.  The determinant of each simplex basis can be
  computed in polynomial time.  Then $\Delta$ can be chosen as the least
  common multiple of all these determinants.
\end{proof}

\begin{lemma}
  \label{lemma:grid-approximation-integral}
  Let $Q\subset \R^d$ be an integral polytope, i.e., all vertices have integer
  coordinates.  Let $M\in\R$ be such that $Q\subseteq \{\, \ve x\in\R^d : |x_i| \leq M \text{ for
  } i=1,\dots,d \,\}$.  Let $\ve x^*\in Q$ and let $\delta>0$.  Then every lattice
  $\frac1k\Z^d$ for $k \geq \frac1\delta (d+1) M$ contains a lattice
  point $\ve x\in Q\cap \frac1k\Z^d$ with $\|\ve x-\ve x^*\|_\infty \leq\delta$.
\end{lemma} 
\begin{proof}
  By Carath\'eodory's Theorem, 
  there exist $d+1$ vertices $\ve x^0, \dots, \ve x^d\in\Z^d$ of~$Q$ and convex multipliers
  $\lambda_0,\dots,\lambda_d$ such that $\ve x^* = \sum_{i=0}^d \lambda_i \ve
  x^i$.  Let $\lambda'_i := \frac1k \lfloor k\lambda_i \rfloor \geq 0$ for
  $i=1,\dots,d$ and $\lambda'_0 := 1 - \sum_{i=1}^d \lambda_i' \geq0$.  
  Then $\ve x := \sum_{i=0}^d \lambda'_i \ve x^i \in Q\cap\frac1k\Z^d$, and
  we have
  \begin{align*}
    \|\ve x-\ve x^*\|_\infty & \leq \sum_{i=0}^d (\lambda'_i - \lambda_i) \| \ve x^i\|_\infty
    \leq (d+1) \frac 1k M 
    \leq \delta.
  \end{align*}
\end{proof}

\section{Bounding techniques for polynomial functions}
\label{section:constant}

Using the results of \autoref{section:qualityofgrid} we are now able to approximate the
mixed-integer optimal point by a point of a suitably fine lattice.  The
question arises how we can use the geometric distance of these two points to
estimate the difference in objective function values.  We prove \autoref{lemma:lipschitz}
that provides us with a local Lipschitz constant for the polynomial to be
maximized.

\begin{lemma}[Local Lipschitz constant]\label{lemma:lipschitz}
  Let $f$ be a polynomial in $d$ variables with maximum
  total degree~$D$. 
  Let $C$ denote the largest absolute value of a coefficient of~$f$.
  Then there exists a Lipschitz constant $L$ such that $|f(\ve x) - f(\ve y)|
  \leq L \|\ve x-\ve y\|_\infty$
  for all $|x_i|, |y_i|\leq M$.  The constant $L$ is $O(D^{d+1} CM^D)$.
\end{lemma}

\begin{proof} 
  Using the usual multi-index notation, let 
  \begin{math}
    f(\ve x) = \sum_{\vealpha\in\mathcal D} c_{\vealpha} \ve x^\vealpha.
  \end{math}
  Let $r = |\mathcal D|$ be the number of monomials of~$f$.
  Then we have
  \begin{displaymath}
    |f(\ve x) - f(\ve y)| \leq \sum_{\vealpha\neq\ve0} |c_\vealpha|\, |\ve x^\vealpha - \ve y^\vealpha|.
  \end{displaymath}
  We estimate all summands separately. 
  Let $\vealpha\neq\ve0$ be an exponent vector with $n := 
  \sum_{i=1}^d \alpha_i \leq D$.  
  Let 
  \begin{displaymath}
    \vealpha = \vealpha^0 \geq \vealpha^1 \geq \dots \geq \vealpha^n = \ve 0
  \end{displaymath}
  be a decreasing chain of exponent vectors with $\vealpha^{i-1} - \vealpha^{i} = \ve
  e^{j_i}$ for $i=1,\dots,n$.  Let $\vebeta^i := \vealpha - \vealpha^i$ for $i=0,\dots,n$. 
  Then $\ve x^\vealpha-\ve y^\vealpha$ can be
  expressed as the ``telescope sum'' 
  \begin{align*}
    \ve x^\vealpha-\ve y^\vealpha
    &= \ve x^{\vealpha^0}\ve y^{\vebeta^0} - \ve x^{\vealpha^1} \ve y^{\vebeta^1} 
    + \ve x^{\vealpha^1} \ve y^{\vebeta^1} - \ve x^{\vealpha^2} \ve y^{\vebeta^2} 
    +- \cdots
    - \ve x^{\vealpha^n} \ve y^{\vebeta^n} \\
    &= \sum_{i=1}^n \left( \ve x^{\vealpha^{i-1}} \ve y^{\vebeta^{i-1}} 
      - \ve x^{\vealpha^{i}} \ve y^{\vebeta^{i}} \right) \\
    &= \sum_{i=1}^n \left( (x_{j_i} - y_{j_i}) \ve x^{\vealpha^{i}} \ve y^{\vebeta^{i-1}} \right). 
  \end{align*}
  Since $ \bigl| \ve x^{\vealpha^{i}} \ve
  y^{\vebeta^{i-1}} \bigr| \leq M^{n-1}$ and $n\leq D$, we obtain 
  \begin{align*}
    | \ve x^\vealpha - \ve y^\vealpha | & \leq D \cdot \|\ve x-\ve y\|_\infty
    \cdot M^{n-1} ,
  \end{align*}
  thus
  \begin{align*}
    |f(\ve x) - f(\ve y)| \leq C r D M^{D-1} \|\ve x-\ve y\|_\infty.
  \end{align*}
  Let $L := C r D M^{D-1} $.  Since $r = O(D^d)$, we have $L = O(D^{d+1} CM^D)$.
\end{proof}

Moreover, in order to obtain an {\small FPTAS}, we need to put differences of
function values in relation to the maximum function value.  
To do this, we need to deal with the special case of polynomials that are constant on
the feasible region; here trivially every feasible solution is optimal.
For non-constant polynomials, we can prove a lower bound on the maximum
function value.  The technique is to bound the difference of the minimum and
the maximum function value on the mixed-integer set from below; if the
polynomial is non-constant, this implies, for a non-negative polynomial, a lower
bound on the maximum function value. We will need a simple fact about the
roots of multivariate polynomials. 

\begin{lemma}\label{lemma:toomanyroots}
  Let $f\in\Q[x_1,\dots,x_{d}]$ be a polynomial and let $D$ be the largest
  power of any variable that appears in~$f$.  Then $f=0$ if and only if $f$
  vanishes on the set $\{0,\dots, D\}^d$. 
\end{lemma}
\begin{proof}
  This is a simple consequence of the Fundamental Theorem of Algebra.  See,
  for instance, \cite[Chapter 1, \S1, Exercise 6\,b]{CoxLittleOShea92}.
\end{proof}

\begin{lemma}\label{lemma:constantness}
  Let $f\in\Q[x_1,\dots,x_{d}]$ be a polynomial with maximum total degree~$D$.
  Let $Q \subset \R^d$ be an integral 
  polytope of dimension $d'\leq d$.  
  Let $k \geq D\, d'$.
  Then $f$ is constant on $Q$ if and only if 
  $f$ is constant on $Q\cap \frac1k \Z^d$.  
\end{lemma}
\begin{proof}
  Let $\ve x^0\in Q\cap\Z^d$ be an arbitrary vertex of~$Q$.  There exist
  vertices $\ve x^1, \dots, \ve x^{d'}\in Q\cap\Z^d$ such that the vectors
  $\ve x^1-\ve x^0, \dots, \ve x^{d'}-\ve x^0\in\Z^d$ are linearly
  independent.  By convexity, $Q$ contains the parallelepiped 
  \begin{displaymath}
    S := 
    \left\{\, \ve x^0 + \textstyle\sum_{i=1}^{d'} \lambda_i (\ve x^i-\ve x^0) :
      \text{$\lambda_i \in [0,\tfrac1{d'}]$ for $i=1,\dots, d'$} \,\right\}.
  \end{displaymath}
  We consider the set 
  \begin{align*}
    S_k &= \tfrac1k \Z^d \cap S
    \supseteq 
    \left\{\, \ve x^0 + \textstyle\sum_{i=1}^{d'}
      \frac{n_i}k (\ve x^i-\ve x^0) : \text{$n_i \in \{0,1,\dots,D\}$ for $i=1,\dots, d'$}
      \,\right\}. 
  \end{align*}
  Now if there exists a $c\in\R$ with $f(\ve x) = c$ for $\ve x\in Q\cap
  \frac1k \Z^d$, then all the points in $S_k$ are roots of the polynomial 
  $f - c$, which has only maximum total degree $D$.  
  By \autoref{lemma:toomanyroots} (after an affine transformation), $f-c$ is
  zero on the affine hull of $S_k$; 
  hence $f$ is constant on $Q$. 
\end{proof}

\begin{theorem}\label{lemma:mixed-constantness}
  Let $f\in\Z[x_1,\dots,x_{d_1},z_1,\dots,z_{d_2}]$.  
  Let $P$ be a rational convex polytope, and let $\Delta$
  be the number from \autoref{lemma:integral-scaling}.
  Let $m = k\Delta$ with $k \geq D\, d_1$, $k\in\Z$.
  Then $f$ is constant on the feasible region $P\cap \bigl(\R^{d_1}\times \Z^{d_2}\bigr)$ if and
  only if $f$ is constant on $P\cap \bigl( \frac1m \Z^{d_1} \times \Z^{d_2}
  \bigr)$.  If $f$ is not constant, then 
  \begin{equation}\label{eq:lower-bound-for-range}
    \bigl| f(\ve x_{\max}, \ve z_{\max})
    - f(\ve x_{\min},\ve z_{\min}) \bigr| \geq m^{-D},
  \end{equation}
  where $(\ve x_{\max}, \ve z_{\max})$ is an optimal solution to the
  maximization problem over the feasible region $P\cap \bigl(\R^{d_1}\times
  \Z^{d_2}\bigr)$ and $(\ve x_{\min}, \ve 
  z_{\min})$ is an optimal solution to the minimization problem.
\end{theorem}
\begin{proof}
  Let $f$ be constant on $P\cap \bigl( \frac1m \Z^{d_1} \times \Z^{d_2}
  \bigr)$.  For fixed integer part $\ve z\in\Z^{d_2}$, we consider the
  polytope $\Delta P_{\ve z} = \bigl\{\, \Delta\ve x :
  (\ve x,\ve z) \in P\,\bigr\}$, which is a slice of $P$ scaled to become an
  integral polytope.  By applying \autoref{lemma:constantness} with $k = (D+1)
  d$ on every polytope $\Delta P_{\ve z}$, we obtain that $f$ is constant on every slice~$P_{\ve z}$. 
  Because $f$ is also constant on the set $P\cap \bigl( \frac1m \Z^{d_1} \times \Z^{d_2}
  \bigr)$, which contains a point of every non-empty slice $P_{\ve z}$,
  it follows that $f$ is constant on $P$. 
  
  If $f$ is not constant, there exist $(\ve x^1, \ve z^1)$, $(\ve x^2, \ve
  z^2) \in P\cap\bigl( \frac1m \Z^{d_1} \times \Z^{d_2} \bigr)$ with $f(\ve x^1,\ve z^1)
  \neq f(\ve x^2,\ve z^2)$.  
  By the integrality of all coefficients of $f$, we
  obtain the estimate
  \begin{align*}
    |f(\ve x^1,\ve z^1) - f(\ve x^2,\ve z^2)| \geq m^{-D}.
  \end{align*}
  Because $(\ve x^1, \ve z^1)$, $(\ve x^2, \ve z^2)$ are both feasible solutions
  to the maximization problem and the minimization problem, this
  implies~\eqref{eq:lower-bound-for-range}. 
\end{proof}

\section{Proof of \autoref{th:mipo-fptas}} \label{section:mainthm}
 
Now we are in the position to prove the main result.

\begin{proof}[Proof of \autoref{th:mipo-fptas}]
  \emph{Part (a).}
  Let $(\ve x^*, \ve z^*)$ denote an optimal solution to the mixed-integer
  problem~\eqref{eq:mixed-problem}.  Let $\epsilon>0$.  We show that, in time
  polynomial in the input length, the maximum total degree,
  and~$\frac1\epsilon$, 
  we can compute a point 
  $(\ve x, \ve z)$ that satisfies
  \mbox{(\ref{eq:mixed-problem-a}--\ref{eq:mixed-problem-b})}
  such that 
  \begin{equation}
    \label{eq:epsilon-approximation}
    |f(\ve x, \ve z) - f(\ve x^*, \ve z^*)| \leq \epsilon f(\ve x^*, \ve z^*).
  \end{equation}

  First we note that we can restrict ourselves to the case of polynomials
  with integer coefficients, simply by multiplying $f$ with the least common
  multiple of all denominators of the coefficients.
  We next establish a lower bound on $f(\ve x^*, \ve z^*)$. 
  To this end, let $\Delta$ be the integer from
  \autoref{lemma:integral-scaling}, which can be computed in
  polynomial time.  By \autoref{lemma:mixed-constantness} with $m = D\, d_1
  \Delta$, either $f$ is constant on the feasible region, or 
  \begin{equation}\label{eq:lower-bound-for-optimum}
    f(\ve x^*, \ve z^*) \geq (D\, d_1 \Delta)^{-D},
  \end{equation}
  where $D$ is the maximum total degree of~$f$.
  Now let 
  \begin{equation}
    \delta:=\frac{\epsilon}{2(D d_1 \Delta)^D L(C,D,M)}
  \end{equation}
  and let
  \begin{equation}
    m := \Delta \left\lceil \frac2\epsilon (D d_1 \Delta)^D  L(C,D,M) (d_1+1) M
    \right\rceil, 
  \end{equation}
  where $L(C,D,M)$ is the Lipschitz constant from \autoref{lemma:lipschitz}.
  Then we have $m \geq \Delta \frac1\delta (d_1+1) M$, so 
  by \autoref{lemma:grid-approximation-lowerdim}, there
  is a point $(\round[\delta]{\ve x^*}, \ve z^*) \in P \cap \bigl( \frac1m
  \Z^{d_1}\times \Z^{d_2}\bigr)$  
  with $\bigl\| \round[\delta]{\ve x^*} - \ve x^* \bigr\|_\infty \leq \delta$. 
  Let $(\ve x^m, \ve z^m)$ denote an optimal solution to the grid problem~\eqref{eq:grid-problem}.  
  Because $(\round[\delta]{\ve x^*}, \ve z^*)$ is a feasible solution to the grid problem~\eqref{eq:grid-problem}, 
  we have 
  \begin{equation}
    f(\round[\delta]{\ve x^*}, \ve z^*) \leq f(\ve x^m, \ve z^m) \leq f(\ve x^*, \ve z^*).
  \end{equation}
  Now we can estimate
  \begin{align}
    \bigl| f(\ve x^*, \ve z^*) - f(\ve x^m, \ve z^m) \bigr| 
    & \leq \bigl| f(\ve x^*, \ve z^*) - f(\round[\delta]{\ve x^*}, \ve z^*) \bigr| \nonumber\\
    & \leq L(C,D,M) \, \bigl\| \ve x^* - \round[\delta]{\ve x^*} \bigr\|_\infty \nonumber\\
    & \leq L(C,D,M) \, \delta \nonumber\\
    & = \frac\epsilon2 {(D\, d_1 \Delta)}^{-D} \nonumber\\
    & \leq \frac\epsilon2 f(\ve x^*, \ve z^*), \label{eq:estimate-gridopt-to-optimum}
  \end{align}
  where the last estimate is given by~\eqref{eq:lower-bound-for-optimum} in
  the case that $f$ is not constant on the feasible region.  On the other
  hand, if $f$ is constant, the
  estimate~\eqref{eq:estimate-gridopt-to-optimum} holds trivially.

  By \autoref{cor:gridproblem-approximation}
  we can compute a point $(\ve 
  x^m_{\epsilon/2}, \ve z^m_{\epsilon/2})\in P\cap\bigl(\frac1m\Z^{d_1}\times
  \Z^{d_2}\bigr)$ such
  that 
  \begin{equation}\label{eq:estimate-gridopt-to-gridapprox}
    (1-\tfrac\epsilon2) f(\ve x^m,\ve z^m) \leq  f(\ve x^m_{\epsilon/2}, \ve
    z^m_{\epsilon/2}) \leq f(\ve x^m, \ve z^m)
  \end{equation}
  in time polynomial in $\log m$, the encoding 
  length of $f$ and~$P$, the maximum total degree~$D$, and $1/\epsilon$. 
  Here $\log m$ is bounded by a polynomial in $\log M$, $D$ and $\log C$,
  so we can compute $(\ve x^m_{\epsilon/2}, \ve z^m_{\epsilon/2})$ 
  in time polynomial in the input
  size, the maximum total degree~$D$, and $1/\epsilon$.
  Now we can estimate, using \eqref{eq:estimate-gridopt-to-gridapprox}
  and~\eqref{eq:estimate-gridopt-to-optimum}, 
  \begin{align*}
    f(\ve x^*, \ve z^*) - f(\ve x^m_{\epsilon/2}, \ve z^m_{\epsilon/2}) 
    & \leq f(\ve x^*,\ve z^*) - (1-\tfrac\epsilon2) f(\ve x^m,\ve z^m) \\
    & = \tfrac\epsilon2 f(\ve x^*,\ve z^*) + (1-\tfrac\epsilon2) \bigl( f(\ve
    x^*,\ve z^*) - f(\ve x^m, \ve z^m)\bigr) \\
    & \leq \tfrac\epsilon2 f(\ve x^*,\ve z^*) + \tfrac\epsilon2 f(\ve x^*,\ve z^*) \\
    & = \epsilon f(\ve x^*,\ve z^*). 
  \end{align*}
  Hence $f(\ve x^m_{\epsilon/2}, \ve z^m_{\epsilon/2}) \geq (1-\epsilon) f(\ve
  x^*, \ve z^*)$.\medbreak

  \noindent\emph{Part (b).}
  Let the dimension $d\geq2$ be fixed.  We prove that there does not exist a
  {\small PTAS} for the maximization of arbitrary polynomials over
  mixed-integer sets of polytopes. 
  We use the NP-complete problem AN1 on page 249 of \cite{GarJohn79}. This is
  to decide whether, given three positive integers $a,b,c$, there exists a
  positive integer $x<c$ such that $x^2 \equiv a\pmod b$.  This problem is
  equivalent to asking whether the maximum of the quartic polynomial
  function $f(x,y) = - (x^2-a-by)^2$ over the lattice points of the rectangle
  \begin{displaymath}
    P = \biggl\{\,(x,y) :  1 \leq x \leq c-1,\ \frac{1-a}{b} \leq y \leq
    \frac{(c-1)^2-a}{b} \,\biggr\}
  \end{displaymath}
  is zero or not.  If there existed a {\small PTAS} for the
  maximization of arbitrary polynomials over mixed-integer sets of polytopes,
  we could, for any fixed $0<\epsilon<1$, compute in polynomial time a solution
  $(x_\epsilon,y_\epsilon)\in P\cap\Z^2$ with $\bigl|f(x_\epsilon,y_\epsilon)
  - f(x^*,y^*)\bigr| \leq \epsilon \bigl|f(x^*,y^*)\bigr|$, where $(x^*,y^*)$ denotes an
  optimal solution.  Thus, we have $f(x_\epsilon,y_\epsilon)=0$ if and only if
  $f(x^*,y^*)=0$; this means we could solve the problem~AN1 in polynomial time.
\end{proof}

\clearpage
\appendix
\section{Appendix: An FPTAS for the integer case}
\label{section:integer-ftpas}

The first fully polynomial-time approximation scheme for the integer case
appeared in our paper~\cite{deloera-hemmecke-koeppe-weismantel:intpoly-fixeddim}. 
It is based on Alexander Barvinok's theory for encoding all the lattice points
of a polyhedron in terms of short rational functions 
\cite{bar,BarviPom}.  The set $P\cap\Z^d$ is represented
by a Laurent polynomial $g_P(\ve z)=\sum_{\vealpha \in P\cap\Z^d}
\ve z^\vealpha.$ From Barvinok's theory this exponentially-large sum of
monomials $g_P(\ve z)$ can be written as a polynomial-size sum of
rational functions (assuming the dimension $d$ is fixed) of the form:
\begin{equation} \label{eq:aa}
g_P(\ve z) = \sum_{i\in I} {E_i \frac{\ve z^{\ve u_i}} {\prod_{j=1}^d
(1-\ve z^{\ve v_{ij}})}},
\end{equation}
where $I$ is a polynomial-size indexing set, and where
$E_i\in\{1,-1\}$ and $\ve u_i, \ve v_{ij} \in\Z^d$ for all $i$ and $j$.
There is a polynomial-time algorithm for computing this representation
\cite{bar,BarviPom,latte2,latte1}.

By symbolically applying differential operators to the
representation~\eqref{eq:aa}, we can compute a short rational function
representation of the Laurent polynomial 
\begin{equation}
  g_{P,f}(\ve z) = \sum_{\vealpha\in P \cap \Z^d} f(\vealpha) \ve z^{\vealpha}.
\end{equation}
In fixed dimension, 
the size of the expressions occuring in the symbolic calculation can be
bounded polynomially:

\begin{lemma}[Lemma 3.1 of
  \cite{deloera-hemmecke-koeppe-weismantel:intpoly-fixeddim}] 
  \label{operators} Let the dimension $d$ be fixed. 
  Let $g_P(\ve z)=\sum_{\vealpha\in P \cap
    \Z^d}\ve z^{\vealpha}$ be the Barvinok 
  representation of the generating function of $P\cap\Z^d$. Let
  $f\in\Z[x_1,\dots,x_d]$ be a polynomial of maximum total
  degree~$D$. We can compute, in time polynomial in $D$ and the input size, a
  Barvinok representation $g_{P,f}(\ve z)$ 
  for the generating function $\sum_{\vealpha\in P \cap \Z^d} f(\vealpha) \ve z^{\vealpha}.$
\end{lemma}

Now we present the algorithm to obtain bounds $U_k,L_k$ that reach the optimum.  
We make use of the elementary fact that, for a set
$S=\{s_1,\dots,s_r\}$ of non-negative real numbers, 
\begin{equation}
  \max\{s_i : s_i \in S \} = \lim_{k \rightarrow \infty} \sqrt[k]{\textstyle\sum_{j=1}^r s_j^k}.
\end{equation}

\begin{algorithm}[Computation of bounds for the maximization problem]~\smallskip
\label{Algorithm}

\noindent {\em Input:} A rational convex polytope $P \subset \R^d$, 
a polynomial objective $f \in \Z[x_1,\dots,x_d]$ of maximum total degree $D$ 
that is non-negative over $P\cap\Z^d$.\smallskip

\noindent {\em Output:} An increasing sequence of lower bounds $L_k$,
and a decreasing sequence of upper bounds $U_k$ reaching the maximal
function value $f^*$ of $f$ over $P\cap\Z^d$.

\begin{enumerate}
\item  Compute a short rational function expression for
  the generating function $g_P(\ve z)=\sum_{\vealpha\in P\cap\Z^d} \ve
  z^{\vealpha}$.  Using residue techniques, compute $|P \cap
  \Z^d|=g_P(\ve 1)$ from $g_P(\ve z)$.

\item From the rational function $ g_P(\ve z)$ 
  compute the rational function representation of $g_{P,f^k}(\ve z)$ of
  $\sum_{\vealpha\in P \cap \Z^d} f^k(\vealpha) \ve z^\vealpha$ in by 
  \autoref{operators}. Using residue techniques, compute 

\[
L_k:=\sqrt[k]{g_{P,f^k}(\ve 1)/g_{P,f^0}(\ve1)}\;\;\;\text{and}\;\;\;
U_k:=\sqrt[k]{g_{P,f^k}(\ve 1)}. 
\] 
\end{enumerate} 
\end{algorithm}

\begin{theorem}[{\small FPTAS}, Lemma 3.3 and Theorem 1.1
  of~\cite{deloera-hemmecke-koeppe-weismantel:intpoly-fixeddim}]\label{lemma:bounds} 
  Let the dimension $d$ be fixed.  Let $P\subset\R^d$ be a rational convex polytope.
  Let $f$ be a polynomial with integer coefficients and maximum total
  degree~$D$ that is non-negative on $P\cap\Z^d$. 
  \begin{enumerate}[\rm(i)]
\item \autoref{Algorithm} computes the bounds $L_k$, $U_k$ in time polynomial in
  $k$, the input size of $P$ and $f$, and the total degree~$D$. The bounds
  satisfy the following inequality:
$$ 
U_k-L_k \leq f^* \cdot \left(\sqrt[k]{|P \cap \Z^d|}-1 \right).
$$
\item For $k=(1+1/\epsilon)\log({|P \cap \Z^d|})$ (a number bounded by a
  polynomial in the input size),
  $L_k$ is a $(1-\epsilon)$-approximation to the optimal value $f^*$ and it
  can be computed in time polynomial in the input size, the total
  degree~$D$, and $1/\epsilon$. Similarly, $U_k$ gives a
  $(1+\epsilon)$-approximation to $f^*$. 

\item With the same complexity,  by iterated bisection of $P$, we can also find
  a feasible solution $\ve x_\epsilon\in P\cap\Z^d$ with 
  \begin{displaymath}
    \bigl|f(\ve x_\epsilon) - f^*\bigr| \leq \epsilon f^*.
  \end{displaymath}
\end{enumerate}
\end{theorem}

\clearpage

\paragraph{Acknowledgment.}
The first author gratefully acknowledges support from NSF grant DMS-0309694, a
2003 UC-Davis Chancellor's fellow award, the Alexander von Humboldt
foundation, and IMO-Magdeburg.  
The remaining authors were supported
by the European TMR network ADONET 504438.

{\samepage
\bibliographystyle{amsplain}
\bibliography{ipo,iba-bib,weismantel}
}
\end{document}